\documentclass[a4paper,10pt]{article}
\usepackage{graphicx}
\usepackage{fancybox}
\usepackage{amssymb}
\usepackage[all]{xy}
\usepackage{color}

\begin{document}
%
\title{Cross-Entropy method: convergence issues for extended implementation}
%
\author{Fr\'ed\'eric Dambreville\\
D\'el\'egation G\'en\'erale pour l'Armement, DGA/CTA/DT/GIP\\
16 Bis, Avenue Prieur de la C\^ote d'Or\\
Arcueil, F 94114, France\\
Web: {\tt http://www.FredericDambreville.com}\\
Email: {\tt http://email.FredericDambreville.com}}
\maketitle
\begin{abstract}
The cross-entropy method (CE) developed by R. Rubinstein is an elegant practical principle for simulating rare events.
The method approximates the probability of the rare event by means of a family of probabilistic models.
The method has been extended to optimization, by considering an optimal event as a rare event.
CE works rather good when dealing with deterministic function optimization.
Now, it appears that two conditions are needed for a good convergence of the method.
First, it is necessary to have a family of models sufficiently flexible for discriminating the optimal events.
Indirectly, it appears also that the function to be optimized should be deterministic.
The purpose of this paper is to consider the case of partially discriminating model family, and of stochastic functions.
It will be shown on simple examples that the CE could fail when relaxing these hypotheses.
Alternative improvements of the CE method are investigated and compared on random examples in order to handle this issue.
\end{abstract}
\section{Introduction}
The Cross-Entropy method has been developed by R. Rubinstein for the simulation of rare events\cite{rubinstein:book}.
The algorithm iteratively builds a near-optimal importance sampling of the rare event, based on a family of parameterized sampling laws.
The construction of the importance sampling is obtained by iteratively:
\begin{itemize}
\item tossing samples,
\item selecting the samples which are approximating the rare events,
\item relearning the parameters of the sampling law by minimizing its Kulback-Leiber distance (cross-entropy) with the selection,
\item computing the importance weightings.
\end{itemize}
By considering the optimal events related to an objective as rare events, the method has been extended to optimization problems.
\\[5pt]
The cross-entropy method has been implemented successfully on many combinatorial problems.
However, attempted proofs of the method make some assumptions as preliminary requests\cite{rubinstein:book,rubinsteinProof}.
First, the proof has been made in a deterministic context.
Secondly, the closure of the simulation law family should contain the dirac on the optimum (or laws with support on the optimums).
\\[5pt]
The first condition cannot be fulfilled properly, in case of stochastic problem.
The second condition is an obvious requirement.
But there are some cases, where it is not possible to handle all the solutions precisely by the law family.
Indeed, the solutions may not be countable practically; this is typically the case for some dynamic problems (for example, the strategy tree against a deterministic computer chess player).
Both difficulties are encountered in optimal planning with partial observation.
The purpose of this paper is to point out on simple examples, that these hypotheses are necessary for the convergence of the classical CE method.
The questions are:
\begin{itemize}
\item \emph{Does the \emph{classical} CE algorithm solve stochastic problems properly?}
It appears that the quantile selection within the CE may not work properly, without a rather good estimation of the objective functional expectation.
Nevertheless, smoother selection criteria seem to be a possible answer to these difficulties.
\item Assume that the law family closure does not contain all the deterministic solutions.
The CE algorithm will converge to a stochastic approximation of the optimal solution.
\emph{Is this approximation the best possible within the law family?}
Our answer to this question is not absolutely negative.
But it appears that some extensions of the CE, quite usually implemented, will fail on this question.

\end{itemize}
This paper presents some counterexamples to these questions.
In the case of stochastic optimization, tests are done on simple random examples in order to compare the convergence of various CE methods with the global optimum.
\\[5pt]
Next section introduces shortly the principle of the CE method.
Section~\ref{RESIM:DMB:Cont1:Chap3} will consider the case, where the optimal solution is not caught properly by the sampling family.
A counterexample is proposed and studied.
In section~\ref{RESIM:DMB:Cont1:Chap4}, stochastic problems are considered.
Two simple counterexamples are investigated, thus enlightening some typical convergence difficulties.
Different evolutions of the cross-entropy are then compared to the basical method, by generating several random examples.
In particular, a method with smooth sample selection is proposed as a possible alternative for the stochastic problems.
Section~\ref{RESIM:DMB:Cont1:Chap5} concludes.
\section{Basis of the cross-entropy method}
\label{RESIM:DMB:Cont1:Chap2}
The reader interested in CE methods should refer to the tutorial \cite{IEEEisda:boer} and the book \cite{rubinstein:book} on the CE method.
CE algorithms were first dedicated to estimating the probability of rare events.
A slight change of the basic algorithm made it also good for optimization.
We will not focus on the cross-entropy method for simulation, although this primary aspect of the method is quite interesting.
Rather, the CE method for optimization is now presented and discussed.
While there are different evolutions of the primary method related to the choice of the selective rate or to a smooth update, this presentation is restricted to the basical CE method.
By the way, it is not difficult to attest that the counterexamples proposed in sections~\ref{RESIM:DMB:Cont1:Chap3} and~\ref{RESIM:DMB:Cont1:Chap4} still work with these evolutions.
\subsection{General CE algorithm for the optimization}
The Cross Entropy algorithm repeats until convergence the three successive phases in order to maximize a given reward criterion:
\begin{enumerate}
\item Generate samples of random data according to a parameterized random mechanism,
\item Select the best samples according to the reward criterion,
\item Update the parameters of the random mechanism, on the basis of the selected samples.
\end{enumerate}
In the particular case of CE, the update in phase 3 is obtained by minimizing the Kullback-Leibler distance, or cross entropy, between the updated random mechanism and the selected samples.
The next paragraphs describe on a theoretical example how such method can be used in an optimization problem.
\paragraph{Formalism.}
Let be given a function $x\mapsto f(x)$; this function is easily computable.
The value $f(x)$ has to be maximized, by optimizing the choice of $x\in X$.
The function $f$ will be the reward criterion.
\vspace{5pt}\\
Now let be given a family of probabilistic laws, $P_\sigma|_{\sigma\in\Sigma}$\,, applying on the variable $x$.
The family $P$ is the parameterized random mechanism.
\vspace{5pt}\\
Let $\rho\in\,]0,1[$ be a selective rate.
The CE algorithm for $(x,f,P)$ follows the synopsis\,:
\begin{enumerate}
\item Initialize $\sigma\in\Sigma$\,,
\item \label{XX:step2}Generate $N$ samples $x_n$ according to $P_\sigma$\,,
\item Select the $\rho N$ best samples according to the reward criterion $f$\,,
\item Update $\sigma$ as a minimizer of the cross-entropy with the selected samples:
$$\sigma\in\arg\max_{\sigma\in\Sigma}\sum_{n~{\rm selected}}\ln P_\sigma(x_n)\;,$$
\item Repeat from step \ref{XX:step2} until convergence.
\end{enumerate}
\emph{This algorithm requires $f$ to be easily computable and the sampling of $P_\sigma$ to be fast.}
%
\paragraph{Interpretation.}
The CE algorithm tightens the law $P_\sigma$ around the maximizer of $f$.
Then, when the probabilistic family $P$ is well suited to the maximization of $f$\,, it becomes equivalent to find a maximizer for $f$ or to optimize the parameter $\sigma$ by means of the CE algorithm.
The problem is to find a good family, and convergence parameters.
\paragraph{Extensions.}
\subparagraph{Smooth update.}
The method has been extended by implementing a smooth update of the law.
More precisely, assume the set $\{P_\sigma / \sigma\in\Sigma\}$ to be convex, and let $\alpha\in[0,1[$ be a smoothing rate.
The algorithm follows the synopsis\,:
\begin{enumerate}
\item Initialize $\sigma\in\Sigma$\,,
\item \label{XX:2:step2}Generate $N$ samples $x_n$ according to $P_\sigma$\,,
\item Select the $\rho N$ best samples according to the reward criterion $f$\,,
\item Define $\sigma_1$ as a minimizer of the cross-entropy with the selected samples:
$$\sigma_1\in\arg\max_{\sigma_1\in\Sigma}\sum_{n~{\rm selected}}\ln P_{\sigma_1}(x_n)\;,$$
\item Define $\sigma_2$ such that $P_{\sigma_2}=\alpha P_{\sigma}+(1-\alpha)P_{\sigma_1}$\,,
and update $\sigma$ by setting $\sigma:=\sigma_2$\,,
\item Repeat from step \ref{XX:2:step2} until convergence.
\end{enumerate}
\subparagraph{Adaptive parameters.}
The principle is to make the parameters $\alpha$ and $\rho$ dependent of the iteration time of the algorithm or on other contextual informations.
Adaptive parameters appears as a main ingredient in the different proofs of convergence of the method.
\subparagraph{Sampling with rejection.}
In some examples (particularly the \emph{salesman}) considered in the CE tutorial~\cite{IEEEisda:boer}, the laws family $P_\sigma|_{\sigma\in\Sigma}$ does not match the set $X$ of valid values for the variable $x$.
More precisely, there is a set $Y\supsetneq X$ such that $P_\sigma\in\mathcal{P}(Y)$, \emph{i.e.} $P_\sigma$ is defined as a probability over $Y$.
The implementation of such a law family in the CE methods is possible by rejecting the invalid samples generated by $P_\lambda$.
A slight change is implied in the step 2 of the CE algorithm:
\begin{enumerate}
\item[2.] Repeat the subsequent process for any $n\in\{1,\dots,N\}\,$:
\begin{enumerate}
\item Generate a sample $x\in Y$ according to $P_\lambda$\,,
\item If $x\not\in X$, then repeat from step~(a)\,,
\item	\emph{At this step, $x\in X$\,.} Then, set $x_n=x$\,.
\end{enumerate}
\end{enumerate}
There is \emph{no other change} implied to the algorithm.
In particular, the update step is the same: \emph{the update of $P_\lambda$ is done from the selected values of the subset $X$.}
\\[5pt]
At first sight, \emph{this update of the law is questionable in regards to the rejection.}
Indeed, the law to be learned from the samples is $P_\lambda/\sum_{x\in X}P_\lambda(x)$ and not $P_\lambda$\,.
This induces a different result while minimizing the cross-entropy with the selected samples.
\\[5pt]
However, the rejection could also be derived from a parameter adaptation: the idea is to interpret the invalid samples of $Y\setminus X$ as samples with very bad reward.
Then, the classical CE scheme is recovered by adapting the number of samples $N$ and the parameter of selection $\rho$ in order to reject these invalid samples.
\\[5pt]
This last interpretation makes sense, when the process actually converges to a law with a support included in $X$.
This is the case, for example, when the law converges to a dirac around the optimum.
But otherwise, it will be shown in section~\ref{RESIM:DMB:Cont1:Chap3} that the convergence may be biased.
\paragraph{Convergence.}
Different convergence results have been proposed for the method and its extensions~\cite{rubinsteinProof, rubinstein:book,Margolin:proof}.
The convergence needs a proper tuning of the parameters of the algorithm (selecting rate, smoothing, number of samples).
Essentially, these results have been established for the optimization of deterministic functions.
Another issue is the stability of the optimization process, when the family of law, $P_\sigma|_{\sigma\in\Sigma}$\,, does not necessarily match the optimal value properly.
The questions investigated by this paper are:
\begin{itemize}
\item \emph{Does the \emph{classical} CE algorithm solve stochastic problems properly?}
A negative answer is given subsequently.
An evolution of the CE is proposed in order to deal with this problem.
\item Assume that the law family closure does not contain the dirac, or dirac mixture, around the optimal solutions.
\emph{Does the CE process provide the best approximation possible within the law family?}
A partial negative answer is provided in next section, by producing a counterexample based on a sampling law with reject.
This counterexample does not work in the classical scheme of the CE.
It is not clearly answered in this paper, what should be the conditions in the CE process for guaranteeing such stability of the convergence.
But it is sure that one have to be more careful in the choice and manipulation of the family.
\end{itemize}
\section{When the family of laws does not enclose the optimum}
\label{RESIM:DMB:Cont1:Chap3}
The subsequent example is inspired from a convergence flaw diagnosed within a practical trajectory planning experiment; an experiment achieved by Francis Celeste~\cite{francis:paper}, which is working in our team.
%
%
\paragraph{Problem setting.}
It is assumed that an agent has two possible actions: the action \emph{continue} or the action \emph{end}.
Each time the agent decides to \emph{continue}, it receives the reward $+1$ and the process is continued.
When the agent decides to \emph{end}, it still receives the reward $+1$ but the process is terminated.
Thus, the agent has to choose a sequence of action, which is a repetition of the action \emph{continue} terminated by the action \emph{end}:
$$
\mathrm{continue};\mathrm{continue}\dots\mathrm{continue};\mathrm{end}\;.
$$
The reward for a whole sequence of action is $t$, the length of the sequence.
Now, a constraint of length is imposed to the actions.
The sequence of action cannot contain more than $T$ actions, so that $t\le T$.
\paragraph{Optimal solution.}
The optimal solution is obvious.
The agent will do as many action as possible.
Its optimal sequence of action is thus:
$$
\underbrace{\mathrm{continue};\cdots;\mathrm{continue}}_{T\times};\mathrm{end}\;.
$$
The problem is actually a triviality.
But we will see that for some laws family, the CE with rejection will fail in finding the optimal law.
\paragraph{Proposal of a laws family, and convergence issue.}
On such a simple example, the best choice is perhaps a law on the length of the process sequence.
But in fact, this kind of problem could be easily generalized so as to involve more than two possible actions (not only \emph{continue} or \emph{end}).
Then, a Markov chain is generally used for these problems.
In the salesman problem, for example, the actions are the choice for a town; the salesman is solved in \cite{rubinstein:book,IEEEisda:boer} by means of a Markov chain with reject.
A method with reject is investigated subsequently.
\\[5pt]
The purpose is to sample a sequence $(d_\theta|1\le\theta\le t)$ where $1\le t\le T$, $d_\theta=\mathrm{continue}$ for $\theta<t$, and $d_t=\mathrm{end}$\,.
This sampling will be done by means of a reject method:
\begin{itemize}
\item Generate a sample without size constraint: $(d_\theta|1\le\theta\le t)$ where $1\le t$, $d_\theta=\mathrm{continue}$ for $\theta<t$, and $d_t=\mathrm{end}$\,,
\item Reject the sample when $t>T$\,.
\end{itemize}
\emph{Sampling a sequence without size constraint.}
The sampling will be generated uniformly and independently for each step, so that the sampling law of the sequence is characterized by the law $p_\lambda$ for sampling a single action:
$$
p_\lambda(d_\theta=\mathrm{continue})=\lambda
\quad\mathrm{and}\quad
p_\lambda(d_\theta=\mathrm{end})=1-\lambda\;.
$$
The whole process takes into account the ending state, so that the sample generation follows the following synopsis:
\begin{enumerate}
\item Set $t=0$\,,
\item \label{Sampling:process:1:step2}Set $t:=t+1$\,,
\item Generate $d_t$ by means of the law $p_\lambda$\,,
\item Repeat from step~\ref{Sampling:process:1:step2}, until $d_t=\mathrm{end}$\,.
\end{enumerate}
As a consequence, the probability of a full sequence $d=(d_\theta|1\le\theta\le t)$ is given by:
$$
P_\lambda(d)=\lambda^{t-1}(1-\lambda)\;.
$$
\emph{Optimal law.}
The optimal law is the one which yields the best gain expectation for the valid trajectories generated by $P_\lambda$\,.
The gain expectation after rejection is given by:
$$
E_{P_\lambda(\cdot|t\le T)}t=
\frac{\sum_{t=1}^T t \lambda^{t-1}(1-\lambda)}{\sum_{t=1}^T \lambda^{t-1}(1-\lambda)}=
\frac{\sum_{t=1}^T t \lambda^{t-1}}{\sum_{t=1}^T \lambda^{t-1}}\;.
$$
This expectation is maximized when $\lambda=1$\,:
$$
\mbox{Within the family, the optimal law is }P_1\;.
$$
Notice that this optimal distribution is not an optimum for the problem.
The family $P_\lambda|0\le\lambda\le 1$ is not sufficiently rich to handle the optimum of the function.
\\[5pt]
It is sometimes not possible to provide a family able to handle the global optimum of the function.
Then, it is often sufficient to find the optimal distribution among the family.
\emph{Is the CE able to provide such optimal distribution among the family?}
Subsequently, it is shown on the example that the CE (with rejection) does not converge to the optimal law $P_1$.
\\[5pt]
\emph{Updating the law.}
Assume $M=\rho N$ samples $(d^n|1\le n\le M)$  being obtained after a sampling process (with reject) and a selection of the best samples.
Denote $t_n$ the ending time of sequence $d^n$ (\emph{beware:} it is not a power operation!).
\\[5pt]
The parameter $\lambda$ for the upcoming loop of the CE algorithm is obtained by maximizing the distance with the selected samples:
$$
\lambda\in\arg\max\frac1M\sum_{n=1}^M \ln P_\lambda(d^n)\;.
$$
Now:
$$\frac1M\sum_{n=1}^M \ln P_\lambda(d^n)=
\frac1M\sum_{n=1}^M \ln\bigl(\lambda^{t_n-1}(1-\lambda)\bigr)=
\left(\left(\frac1M\sum_{n=1}^M t_n\right)-1\right)\ln\lambda+\ln(1-\lambda)\,.
$$
The maximization then results to the relation:
$$
\left(\left(\frac1M\sum_{n=1}^M t_n\right)-1\right)\frac1\lambda-\frac1{1-\lambda}=0\;.
$$
At last, the following update relation is derived:
\begin{equation}
\lambda=1-M\left/\sum_{n=1}^M t_n\right.\;.
\label{Resim2006:eq:1}
\end{equation}
\emph{Convergence issue.}
Equation~\ref{Resim2006:eq:1} and the rejection constraint imply that $\lambda\le1-\frac1T$ after update.
As a consequence, \emph{the CE does not converge to $P_1$, the optimal distribution among the family.}
In fact, it is even proved by considering the CE process that $\lambda<1-\frac1T$.
Let $P_{\lambda\ast}$ be the law obtained after convergence of the CE.
Then:
$$
E_{P_\lambda^\ast(\cdot|t\le T)}t<
\frac{\sum_{t=1}^T t (1-1/T)^{t-1}}{\sum_{t=1}^T (1-1/T)^{t-1}}\;.
$$
\\[5pt]
Let us consider the simple case $T=2$, and compare the expectations:
$$
E_{P_1(\cdot|t\le T)}t=(1+2)/(1+1)=\frac32
\quad\mbox{and}\quad
E_{P_\lambda^\ast(\cdot|t\le T)}t<(1+2\times\frac12)/(1+\frac12)=\frac43\;.
$$
The difference, at least $11\%$, is not negligible.
\paragraph{Convergence in the CE classical scheme.}
As it has been discussed in section~\ref{RESIM:DMB:Cont1:Chap2}, the update of $\lambda$ within the classical scheme will be obtained by minimizing the cross-entropy of the \emph{conditional} law:
$$
P_\lambda^\ast(d|t\le T)=\frac{\lambda^{t-1}(1-\lambda)}{\sum_{\theta=1}^T\lambda^{\theta-1}(1-\lambda)}=\frac{\lambda^{t-1}(1-\lambda)}{1-\lambda^T}\;.
$$
with the selected samples.
Thus, the update is expressed by:
$$
\lambda\in\arg\max\frac1M\sum_{n=1}^M\ln\frac{\lambda^{t_n-1}(1-\lambda)}{1-\lambda^T}\;.
$$
Defining $\bar t=\frac1M\sum_{n=1}^Mt_n$, the optimization reduces to:
$$
\lambda\in\arg\max\frac{\lambda^{\bar t-1}(1-\lambda)}{1-\lambda^T}\;.
$$
The maximum of this function is not necessarily located at $\lambda=1$\,.
For example, when $\bar t=1$, the optimum is obtained for $\lambda=0$\,.
Now, the function to be optimized could be rewritten:
$$
\frac{\lambda^{\bar t-1}(1-\lambda)}{1-\lambda^T}=\frac{\lambda^{\bar t-T}}{\sum_{k=0}^{T-1}\lambda^{-k}}\;.
$$
Then, it is deduced:
\begin{equation}\label{RESIM2006:MaxEnt:eq:1}
\bar t\ge T\quad\Longrightarrow\quad 1\in\arg\max_{0\le\lambda\le1}\frac{\lambda^{\bar t-1}(1-\lambda)}{1-\lambda^T}\;.
\end{equation}
The equation~(\ref{RESIM2006:MaxEnt:eq:1}) has a clear interpretation: when $\lambda>0$ at initialization and the selective rate $\rho$ is sufficiently small, then the CE algorithm (without rejection) converge to the optimal law $P_1$\,.
As a conclusion, our counterexample fails in the classical CE paradigm.
\paragraph{Discussion.} 
The previous example has shown convergence issue of the CE with reject when the laws family cannot reach the optimum of the function.
This counterexample does not work when using a classical CE scheme.
In general, even when the family cannot handle the optimum exactly, the convergence still works rather well in the classical CE paradigm.
Many questions arise however.
In particular, how to evaluate and enhance the stability of the convergence in regards to the discrimination weakness of the laws family?
\section{When the problem is stochastic}
\label{RESIM:DMB:Cont1:Chap4}
In this section, it is discussed about the convergence of the CE in case of stochastic optimization.
Notice that it is still possible to bring such stochastic problems back to deterministic problems by computing the expectation of the objective function.
But generally, this computation is obtained by simulation and is costy.
A reduction of the cost could be obtained by means of the method described in section~\ref{RESIM2006:STOCHMETH:section}.
\\[5pt]
When the variable to be optimized and the stochastic variable of the system are dependent, the expectation will make necessary the use of a functional abstraction of the variable to be optimized (instead of conditional laws).
This is again somewhat costy.
Moreover, the cost reduction method described in section~\ref{RESIM2006:STOCHMETH:section} is no more feasible (when the variable of the system depends on the variable to be optimized). 
\\[5pt]
The purpose of this section is to consider the stochastic optimization by the CE without computing the expectation of the objective.
It is shown on simple examples that there may be a true convergence difficulty of the CE method in such conditions.
\\[5pt]
In the first subsequent example, the value to be optimized is conditioned by another variable which is stochastic.
In other word, the value to be optimized could be considered as a function of the stochastic variable.
Such problems do not appear classically in the CE literacy, but explain clearly some typical difficulties in the convergence.
The second example is unconditioned and more classical.
These examples will be completed by a study of stochastic optimization problems (here, without conditioning), which will be generated randomly.
Alternative solutions to the classical CE are proposed and compared then.
\subsection{Examples}
\subsubsection{Example 1}
Typically, there is an additional difficulty in evaluating the expectation of the objective function, when the variable to be optimized are conditioned by the variable of the system.
For this reason, we will start by considering this kind of example.
\\[5pt]
Let us consider the following stochastic problem:
\begin{equation}
\begin{array}{@{}l@{}}\displaystyle
f_o\in\arg\max_{f:x\mapsto d}\sum_{x}p(x)V(f(x),x)\;,
\vspace{7pt}\\
\mbox{where }
x\in\{0,1\}\,,\ d\in\{0,1\}\,,\ p(0)=p(1)=\frac12\,,\ V(d,x)=2x+d\;,
\vspace{7pt}\\
\mbox{and }f\mbox{ is a mapping from }x\mbox{ to }d\;.
\end{array}
\label{Resim2006:eq:2}
\end{equation}
This problem could be seen from a probabilistic viewpoint:
\begin{equation}
\begin{array}{@{}l@{}}\displaystyle
h_o\in\arg\max_{h}\sum_{d,x}p(x)h(d|x)V(d,x)\;,
\vspace{7pt}\\
\mbox{where }
x\in\{0,1\}\,,\ d\in\{0,1\}\,,\ p(0)=p(1)=\frac12\,, V(d,x)=2x+d\;,
\vspace{7pt}\\
\mbox{and }h(d|x)\mbox{ is a probability of }d\mbox{ conditionally to }x\;.
\end{array}
\label{Resim2006:eq:3}
\end{equation}
We will apply a cross-entropic method in order to solve the optimization~(\ref{Resim2006:eq:3}).
Notice that the method will differ slightly from usually, since we are handling a conditional laws family.
\paragraph{Direct solve.}
The obvious answer to this problem is $h(0|x)=0$ and $h(1|x)=1$; the optimal gain is $2$.
\paragraph{Cross-entropic solve.}
A cross-entropic procedure is proposed here with quantile selection $\rho=10\%$ (no smooth update, for simplicity) in order to solve~(\ref{Resim2006:eq:3}):
\begin{itemize}
\item Initialize $h$ by $h(0|0)=h(1|0)=h(0|1)=h(1|1)=\frac12$,
\item Make 100 samples and evaluate them by $V$,
\item $[\ast]$ Select the $10\%$ best samples, update $h(\cdot|x)$ from the selected samples, when it is possible.\footnote{Leave $h(\cdot|x)$ unchanged when there are no selected samples conditioned by $x$.}
Reiterate from previous step.
\end{itemize}
Since $V(d_1,1)>V(d_2,0)$ for any choice of $d_i$, it comes that samples $(d,0)$ are (almost) never\footnote{Probability is around $10^{-18}$} selected.
As a consequence, $h(\cdot|0)$ is (almost) never updated and stays unchanged.
Thus, a practical convergence will stale to the solution $h(0|0)=h(1|0)=\frac12\,;\;h(0|1)=0\,;\;h(1|1)=1$, which is sub-optimal.
The expected gain is then $\frac74$.
\\[5pt]
This example contains a specific difficulty: we are indeed optimizing the function $x\mapsto f(x)$ by mean of a conditional law.
By the way, one may argue that $[\ast]$ is not a good updating strategy, since the samples should be selected relatively to each condition $x$.
But this is not possible, when there are many possible conditions $x$ (this is often the case).
\subsubsection{Example 2}
\label{stoch:counter:ex:2}
It could be argued about the previous example that the use of a conditional family is not the classical scheme for applying the CE method.
This forthcoming example will be related to a more classical scheme.
\\[4pt]
Now, let us solve the following stochastic problem:
\begin{equation}
\begin{array}{@{}l@{}}\displaystyle
d_o\in\arg\max_{d}\sum_{x}p(x)V(d,x)\;,
\vspace{7pt}\\
\mbox{where }
x\in\{0,1\}\,,\ d\in\{0,1\}\,,\ p(0)=p(1)=\frac12\;,
\vspace{7pt}\\
\mbox{ and }
V(0,0)=2\,,\ 
V(0,1)=-2\,,\ 
V(1,0)=V(1,1)=1\,.
\end{array}
\label{Resim2006:eq:4}
\end{equation}
From a CE viewpoint, the problem becomes:
\begin{equation}
\begin{array}{@{}l@{}}\displaystyle
h_o\in\arg\max_{h}\sum_{x}p(x)h(d)V(d,x)\;,
\vspace{7pt}\\
\mbox{where }
x\in\{0,1\}\,,\ d\in\{0,1\}\,,\ p(0)=p(1)=\frac12\;,
\vspace{7pt}\\
V(0,0)=2\,,\ 
V(0,1)=-2\,,\ 
V(1,0)=V(1,1)=1\,.
\vspace{7pt}\\
\mbox{and }h(d)\mbox{ is a probability of }d\;.
\end{array}
\label{Resim2006:eq:5}
\end{equation}
\paragraph{Direct solve.}
The optimal solution of~(\ref{Resim2006:eq:5}) is of course $h_o(0)=0$ and $h_o(1)=1$\,, resulting in the gain $1$.
\paragraph{Cross-entropic solve.}
A cross-entropic procedure is proposed here with quantile selection $\rho=10\%$ (no smooth update, for simplicity) in order to solve~(\ref{Resim2006:eq:5}):
\begin{itemize}
\item Initialize $h$ by $h(0)=h(1)=\frac12$,
\item Make 100 samples and evaluate them by $V$,
\item Select the $10\%$ best samples, update $h$ from the selected samples.
Reiterate from previous step.
\end{itemize}
Since $V(0,0)>V(d,x)$ for any $(d,x)\ne (0,0)$\,, it comes that the samples $(d,x)\ne (0,0)$ are (almost) never selected.
As a consequence, the selected samples will be $(0,0)$ from the beginning of the CE process.
Consequently, the CE process will converge to the sub-optimal solution $h_\ast(0)=1$ and $h_\ast(1)=0$\,, thus resulting in the gain $0$.
\subsubsection{Discussion.}
The two previous examples are enlightening.
It appears clearly that the selection scheme of the CE (selection of a quantile) does not work properly, in regards to a stochastic objective.
Indeed, some configurations of the problem, which are sampled by the law of the system but not by us, will be automatically discarded by the quantile selection process.
By discarding these cases, a convergence bias is generated.
\subsection{Alternative methods}
\subsubsection{Computing the expectation (reduced cost)}\label{RESIM2006:STOCHMETH:section}
This method is not exactly an alternative: it is costy.
But it will be provided as a reference for the test comparison.
The idea is to replace the stochastic objective function $V(d,x)$ by an estimation of its expectation.
This expectation is obtained by sampling over $x$ according to the law $p$ of the system.
More samples are used, more accurate is the estimation.
Here, we are using the same samples of $x$ for computing the expected gain of the samples $d_n$.
This will reduce greatly the complexity.
But such method is not feasible, when the variables $x$ and $d$ are dependent.
The whole algorithm is explained subsequently:
\begin{enumerate}
\item Initialize $h$\,,
\item \label{XX:step76} Generate $N$ samples $d_n$ according to $h$\,,
\item Generate $K$ samples $x_k$ according to $p$\,,
\item Evaluate each sample $d_n$ by the estimated expectation $v_n=\sum_{k=1}^KV(d_n,x_k)$\,,
\item Select the $\rho N$ best samples $d_n$ according to the expectation $v_n$\,,
\item Update $h$ as a minimizer of the cross-entropy with the selected samples:
$$h\in\arg\max\sum_{n~{\rm selected}}\ln h(d_n)\;,$$
\item Repeat from step \ref{XX:step76} until convergence.
\end{enumerate}
\subsubsection{Using another selection scheme for the CE}
The idea here is to change the selection scheme of the CE.
The stochastic objective function $V(d,x)$ is directly used here.
As in section~\ref{stoch:counter:ex:2}, the stochastic pair $(d,x)$ is sampled and evaluated at the same time.
\paragraph{Selection scheme.}
Assume $N$ samples $(d_n,x_n)$ being evaluated by $v_n=V(d_n,x_n)$\,.
It is defined a non decreasing function $R$, which will characterize the importance $R(v_n)$ of each sample $(d_n,x_n)$.
The update of $h$ will be computed as a maximizer of the cross entropy with the discrete weighted distribution $\left(d_n,\frac{R(v_n)}{\sum_{n=1}^NR(v_n)}\right)$\,.
\paragraph{Algorithm.}
The whole algorithm is explained subsequently:
\begin{enumerate}
\item Initialize $h$\,,
\item \label{XX:step71} Generate $N$ samples $d_n$ according to $h$ and $N$ samples $x_n$ according to $p$\,,
\item Evaluate each sample pair $(d_n,x_n)$ by $v_n=V(d_n,x_n)$\,,
\item Update $h$ as a minimizer of the cross-entropy with the weighted samples:
$$h\in\arg\max\sum_{n=1}^NR(v_n)\ln h(d_n)\;,$$
\item Repeat from step \ref{XX:step71} until convergence.
\end{enumerate}
This selection scheme is called \emph{smooth selection scheme}.
Notice that the quantile selection of Rubinstein is a particular case of the \emph{smooth selection scheme}, where the function $R$ is a heavyside function pointed on the quantile.
\subsection{Method comparison by  means of Randomly generated tests}
The three methods, basic CE; CE with expectation computation; and smooth selection scheme, have been compared on random problems.
The method for creating the problems is simple:
\begin{itemize}
\item There are $100$ possibles states for $d$ and for $x$, that is $d,x\in\{1,\dots,100\}$\,,
\item The parameters $V(d,x)\in]0,1]$ are generated randomly, according to the uniform law, for any $d$ and any $x$,
\item The probability $p$ is generated randomly, according to the uniform law (that is the 99-dimensions vector characterizing $p$ is generated uniformly),
\end{itemize}
Notice that it is quite easy to solve these problems, by enumerating the cases.
\\[5pt]
The test has been executed $1000$ times.
The parameters of the algorithm are:
\begin{itemize}
\item $K=N=100$ and $\rho=10\%$,
\item The update is smoothed by $\alpha=0.9$ (\emph{i.e.} the innovation is $10\%$),
\item The importance function $R$ is defined by $R(v_n)=v_n$\,.
\end{itemize}
The following table gives the percentage of the optimum achieved by each method.
These results are averaged over the $1000$ executed tests, and the variance is given.
$$
\begin{array}{c||c|c|}
\mbox{Optimal percentage}&\mbox{Mean}&\mbox{Variance}
\\\hline
\mbox{Basic CE}&93.9\%&3.7\%
\\\hline
\mbox{Expectation}&99\%&0.4\%
\\\hline
\mbox{Smooth scheme}&99.1\%&0.7\%
\end{array}
$$
The convergence speed of the expectation CE and the smooth selection CE was comparable.
Since the expectation is computed with reduced cost, the methods run with similar computation cost.
\section{Conclusion}
\label{RESIM:DMB:Cont1:Chap5}
This paper has investigated the convergence issues of the cross-entropy method when relaxing the constraints of use.
A counterexample has been found for the CE with reject, when the laws family used for the CE is too weak and does not contain the optimum dirac.
Counterexamples have been found when optimizing a stochastic objective function.
Weakness of the family and stochastic objective are very important context of use of the CE algorithm.
By the way, both difficulties are encountered when optimizing a control with partial observation\cite{dambreville}.
An alternative evolution of the CE has been proposed for the stochastic optimization.
It is based on a smooth scheme for the sample selection.
The convergence of weak laws family is still an unsolved question.
Next works will focus on this difficult problem.
Moreover, the proof of convergence of the  smooth selection scheme will be investigated; at this time, this method has been evaluated only by experimental means.
\end{document}